\documentclass[11pt,reqno]{amsart}
\usepackage{amsmath,amsfonts, amssymb, amsthm}
  \usepackage{graphics} 
  \usepackage{epsfig} 
\usepackage{graphicx}  \usepackage{epstopdf}
 \usepackage[colorlinks=true]{hyperref}
\hypersetup{urlcolor=blue, citecolor=red}

\newtheorem{theorem}{Theorem}[section]

\newtheorem{lemma}[theorem]{Lemma}

\newtheorem{conjecture}{Conjecture}[section]

\theoremstyle{definition}

\newtheorem{remark}{Remark}[section]

\def\pmod #1{\ ({\rm{mod}}\ #1)}
\def\Z{\Bbb Z}
\def\N{\Bbb N}

\def\Q{\Bbb Q}

\def\l{\left}
\def\r{\right}
\def\bg{\bigg}
\def\({\bg(}
\def\){\bg)}
\def\t{\text}
\def\f{\frac}

\def\mo{{\rm{mod}\ }}
\def\pmod#1{\ (\mo\ #1)}

\def\ls{\leqslant}
\def\gs{\geqslant}

\def\sm{\setminus}
\def\bi{\binom}

\def\eq{\equiv}

\def\Proof{\noindent{\it Proof}}

\begin{document}
\hbox{Chin. Quart. J. Math. 40 (2025), no.\,4, 372--392.}
\medskip

\title[New series involving binomial coefficients (III)]
      {New series involving binomial coefficients (III)}
\author[Zhi-Wei Sun]{Zhi-Wei Sun}
\address{School of Mathematics, Nanjing
University, Nanjing 210093, People's Republic of China}
\email{{\tt zwsun@nju.edu.cn}
\newline\indent
{\it Homepage}: {\tt http://maths.nju.edu.cn/\lower0.5ex\hbox{\~{}}zwsun}}

\keywords{Binomial coefficient, combinatorial identity, infinite series, Kronecker symbol, $L$-function.
\newline \indent Supported by the Natural Science Foundation of China (grant no. 12371004).}

\begin{abstract} We evaluate some series with summands involving a single binomial coefficient
$\bi{6k}{3k}$. For example, we prove that
$$\sum_{k=0}^\infty\f{(63k^2+78k+22)8^k}{(2k+1)(6k+1)(6k+5)\bi{6k}{3k}}=\f{3\pi}2.$$
Motivated by Galois theory, we introduce the so-called Duality Principle for irrational series of Ramanujan's type or Zeilberger's type, and apply it to find 26 new irrational series
identities. For example, we conjecture that
\begin{align*}&\sum_{k=1}^\infty\frac{(32(91\sqrt{33}-523))^{k}}{k^3\binom{2k}k^2\bi{3k}k}
\left((91\sqrt{33}+891)k-33\sqrt{33}-225\right)
\\&\qquad=320\left(\frac{11}3\sqrt{33}L_{-11}(2)-27L_{-3}(2)\right),
\end{align*}
where $ L_{d}(2)=\sum_{k=1}^\infty\frac{(\frac{d}k)}{k^2}$
for any integer $d\equiv0,1\pmod4$ with $(\frac{d}k)$ the Kronecker symbol.
\end{abstract}
\maketitle

\section{Introduction}
\setcounter{equation}{0}
 \setcounter{conjecture}{0}
 \setcounter{theorem}{0}
 \setcounter{proposition}{0}

 In 1914, S. Ramanujan \cite{R} proposed $17$ conjectural series for $1/\pi$ including the identity
 $$\sum_{k=0}^\infty(4k+1)\f{\bi{2k}k^3}{(-64)^k}=\f2{\pi}$$
 first found by G. Bauer in 1859.  The classical Ramanujan-type series for $1/\pi$ has the form:
$$\sum_{k=0}^\infty(ak+b)\f{c_k}{m^k}=\f{\sqrt{d}}{\pi},$$
where $a,b,m\in\Z$ with $am\not=0$, $d$ is a positive rational number, and $c_k$ with $k\in\N=\{0,1,2,\ldots\}$ is among
\begin{equation}\label{three}\bi{2k}k^3,\ \bi{2k}k^2\bi{3k}k,\ \bi{2k}k^2\bi{4k}{2k},\ \bi{2k}k\bi{3k}k\bi{6k}{3k}.
\end{equation}
 Such series are called {\it rational series} since $a,b,m$ are rational numbers.
 There are 36 known rational Ramanujan-type series for $1/\pi$ (cf.\cite [Chapter 14]{Co17}), and it is widely believed that the list is complete. By Stirling's formula,
 $$\lim_{k\to\infty}\root k\of{c_k}=\begin{cases}64&\t{if}\ c_k=\bi{2k}k^3,
 \\108&\t{if}\ c_k=\bi{2k}k^2\bi{3k}k,
 \\256&\t{if}\ c_k=\bi{2k}k^2\bi{4k}{2k},
 \\1728&\t{if}\ c_k=\bi{2k}k\bi{3k}k\bi{6k}{3k}.
 \end{cases}$$

 In 1993, D. Zeilberger \cite{Z} used the WZ method to obtain the identity
 \begin{equation}\label{Z}\sum_{k=1}^\infty\f{21k-8}{k^3\bi{2k}k^3}=\f{\pi^2}2.\end{equation}
 Along this line, J. Guillera \cite{G08} deduced the following identities of Zeilberger's type:
  \begin{equation}\label{G1}\sum_{k=1}^\infty\f{(4k-1)(-64)^k}{k^3\bi{2k}k^3}=-16G,\ \sum_{k=1}^\infty\f{(3k-1)(-8)^k}{k^3\bi{2k}k^3}=-2G,
  \end{equation}
 and
 \begin{equation}\label{G2}\sum_{k=1}^\infty\f{(3k-1)16^k}{k^3\bi{2k}k^3}
 =\f{\pi^2}2,\end{equation}
 where $G$ denotes the Catalan constant $\sum_{k=0}^\infty (-1)^k/(2k+1)^2.$
 In 2010, the author \cite{S11} conjectured the following similar identities:
 \begin{equation}\label{S1}\sum_{k=1}^\infty\f{(10k-3)8^k}{k^3\bi{2k}k^2\bi{3k}k}=\f{\pi^2}2,\
 \sum_{k=1}^\infty\f{(11k-3)64^k}{k^3\bi{2k}k^2\bi{3k}k}=8\pi^2, \ \sum_{k=1}^\infty\f{(35k-8)81^k}{k^3\bi{2k}k^2\bi{4k}{2k}}=12\pi^2,
 \end{equation}
 \begin{equation}\label{S2}\sum_{k=1}^\infty\f{(15k-4)(-27)^{k-1}}{k^3\bi{2k}k^2\bi{3k}k}=K
 \ \t{and}\ \sum_{k=1}^\infty\f{(5k-1)(-144)^k}{k^3\bi{2k}k^2\bi{4k}{2k}}=-\f{45}2K,
 \end{equation}
where
$$K:=\sum_{k=1}^\infty\f{(\f k3)}{k^2}=\sum_{n=0}^\infty\l(\f1{(3n+1)^2}-\f1{(3n+2)^2}\r)$$
with $(\f k3)$ the Legendre symbol. All the five conjectural identities have been confirmed, see
\cite{HP,Gu,GR}. In 2024,  F. Calegari, V. Dimitrov and Y. Tang \cite{CDT} proved that $K$ is irrational while the Catalan constant $G$ is not yet known to be irrational. Note that \cite[(5.1)]{SII} gives the fastest way to compute the value of $K$.

In 1974 R. W. Gosper discovered the identity
$$\sum_{k=0}^\infty\f{25k-3}{2^k\bi{3k}k}=\f{\pi}2$$
which was later proved in \cite{AKP} via integration.
 The author \cite{S24} evaluated some series with summands involving a single binomial coefficient $\bi{3k}k$ or $\bi{4k}{2k}$; for example, the author \cite{S24} proved that
 $$\sum_{k=0}^\infty\frac{(49k+1)8^k}{3^k\binom{3k}k}=81+16\sqrt3\,\pi
\ \ \t{and}\ \
\sum_{k=0}^\infty\frac{10k-1}{\binom{4k}{2k}}=\frac{4\sqrt 3}{27}\pi.$$
  In \cite{SII} the author focused on series of the type
 $$\sum_{k=1}^\infty\f{ak^2+bk+c}{k(3k-1)(3k-2)m^k\bi{4k}k},$$
 where $a,b,c$ and $m\not=0$ are rational numbers; for example, the author \cite{SII} established the identities
 $$\sum_{k=1}^\infty\frac{(5k^2-4k+1)8^{k}}{k(3k-1)(3k-2)\binom{4k}k}=\frac{3}2\pi$$
and $$\sum_{k=1}^\infty\frac{415k^2-343k+62}{k(3k-1)(3k-2)(-8)^k\binom{4k}k}=-3\log2.$$

In 2014 W. Chu and W. Zhang \cite{CZ} used hypergeometric transformations to deduce many series identities involving binomial coefficients. In particular, Examples 58, 115 and 117 of \cite{CZ} yield some series identities
with summands involving a single binomial coefficient $\bi{6k}{3k}$. For example, \cite[Example 115]{CZ} gives the identity
\begin{equation}\label{CZ115}\sum_{k=0}^\infty\f{21k^2+27k+8}{(2k+1)(6k+1)(6k+5)\bi{6k}{3k}}
=\f{8\pi}{9\sqrt3}.
\end{equation}

In this paper we deduce the following new results.

\begin{theorem} \label{Th-8^k} We have the identity
 \begin{equation}\label{pi}\sum_{k=0}^\infty\f{(63k^2+78k+22)8^k}{(2k+1)(6k+1)(6k+5)\bi{6k}{3k}}=\f{3\pi}2.
 \end{equation}
 \end{theorem}

 \begin{theorem}\label{Th1.2} We have \begin{align}\label{42k-1}\sum_{k=0}^\infty\f{(6k-1)(42k-1)\bi{6k}{3k}}{4096^k}&=\f{32}{27}\sqrt3,
  \\\label{42k+5}\sum_{k=0}^\infty\f{(42k+5)\bi{6k}{3k}}{(6k-5)4096^k}&=-\f49\sqrt3,
  \end{align}
  and
  \begin{equation}\label{252}\sum_{k=0}^\infty\f{(252k^2-120k+5)\bi{6k}{3k}}{(2k-1)(6k-1)(6k-5)4096^k}=-\f{\sqrt3}2.
  \end{equation}
 \end{theorem}

In the next section, we will prove Theorem \ref{Th-8^k} via integration, and deduce Theorem \ref{Th1.2} by using the identity
 \begin{equation}\label{CZ4096}\sum_{k=0}^\infty\f{k(252k^2-264k+61)\bi{6k}{3k}}{(2k-1)(6k-1)(6k-5)4096^k}
 =\f1{12\sqrt3}\end{equation}
 given by \cite[Example 117]{CZ}.

Recall that the usual harmonic numbers are defined by
$$H_n:=\sum_{0<k\ls n}\f1k\ \ (n=0,1,2,\ldots).$$
For any integer $m\gs2$, the harmonic numbers of order $m$ are given by
$$H_n^{(m)}:=\sum_{0<k\ls n}\f1{k^m}\ \ (n=0,1,2,\ldots).$$
Motivated by \cite{S15,GL,W,WR}, the author \cite{harmonic,S24,SII} posed many conjectural series identities
involving harmonic numbers. For example, inspired by the Ramanujan series
$$\sum_{k=0}^\infty(42k+5)\f{\bi{2k}k^3}{4096^k}=\f{16}{\pi},$$
the author \cite[(110)]{harmonic} conjectured the identity
\begin{equation*} \sum_{k=0}^\infty(42k+5)\f{\bi{2k}k^3}{4096^k}\l(H_{2k}^{(3)}-\f{43}{352}H_k^{(3)}\r)
 =\f{555}{77}\cdot\f{\zeta(3)}{\pi}-\f{32}{11}G
 \end{equation*}
 which remains open. In Section 3, we will pose our conjectures on series involving $\bi{6k}{3k}$
 some of which also involve harmonic numbers.

 In Section 4, we will introduce the so-called Duality Principle for irrational series of Ramanujan's type and Zeilberger's type. Motivated by this Duality Principle, in Sections 5 and 6 we will propose totally 26 irrational conjectural series of Ramanujan's type or Zeilberger's type.

 All the series in this paper converge at geometric rates and hence evaluations of them can be easily checked numerically via {\tt Mathematica}.

 \section{Proofs of Theorems \ref{Th-8^k} and \ref{Th1.2}}
\setcounter{equation}{0}
 \setcounter{conjecture}{0}
 \setcounter{theorem}{0}
 \setcounter{proposition}{0}

 \begin{lemma}\label{Lem2.1} We have
 \begin{equation}\label{3pi/2}\sum_{k=0}^\infty\f{(63k^2+18k+2)8^k}{(6k+1)\bi{6k}{3k}}=\f32\pi+4.
 \end{equation}
 \end{lemma}
 \Proof. For any $a>0$ and $b>0$, it is well known that
 $$B(a,b):=\int_{0}^1x^{a-1}(1-x)^{b-1}dx$$
 coincides with $\Gamma(a)\Gamma(b)/\Gamma(a+b)$. Thus
 \begin{align*}&\ \ \ \sum_{k=0}^\infty\f{(63k^2+18k+2)8^k}{(6k+1)\bi{6k}{3k}}
 \\&=\sum_{k=0}^\infty(63k^2+18k+2)8^k\f{\Gamma(3k+1)^2}{\Gamma(6k+2)}
 \\&=\sum_{k=0}^\infty (63k^2+18k+2)8^k B(3k+1,3k+1)
 \\&=\int_0^1\sum_{k=0}^\infty (63k^2+18k+2)(2x(1-x))^{3k}dx
 \end{align*}

 For $|z|<1$, by calculus we have
 $$\sum_{k=0}^\infty z^k=\f1{1-z},\ \sum_{k=1}^\infty kz^{k-1}=\f d{dz}(1-z)^{-1}=\f 1{(1-z)^2},$$
 and
 $$\sum_{k=0}^\infty k(k-1)z^{k}=z^2\f d{dz}(1-z)^{-2}=\f{2z^2}{(1-z)^3}.$$
 Thus, for $0\ls x\ls1$ we have
\begin{align*}&\ \ \ \sum_{k=0}^\infty (63k^2+18k+2)(2x(1-x))^{3k}
\\&=\sum_{k=0}^\infty(63k(k-1)+81k+2)(2x(1-x))^{3k}
\\&=\f{126(2x(1-x))^6}{(1-(2x(1-x))^3)^3}+\f{81(2x(1-x))^3}{(1-(2x(1-x))^3)^2}+\f2{(1-(2x(1-x))^3)}
\\&=\f{47(2x(1-x))^6+77(2x(1-x))^3+2}{(1-(2x(1-x))^3)^3}
\\&=\f{8(47(1-t^2)^6+616(1-t^2)^3+128)}{(8-(1-t^2)^3)^3},
\end{align*}
where $t=2x-1$. Combining this with the last paragraph, we obtain
\begin{align*}&\ \ \ \sum_{k=0}^\infty\f{(63k^2+18k+2)8^k}{(6k+1)\bi{6k}{3k}}
\\&=8\int_{-1}^1\f{47(1-t^2)^6+616(1-t^2)^3+128}{2(8-(1-t^2)^3)^3}dt
\\&=8\int_{0}^1\f{47(1-t^2)^6+616(1-t^2)^3+128}{(8-(1-t^2)^3)^3}dt.
\end{align*}
For
$$f(t):=\f{48t(t^4-3t^2+10)}{(t^6-3t^4+3t^2+7)^2}+\f{t(3t^4-10t^2-25)}{t^6-3t^4+3t^2+7}+3\arctan t,$$
it is easy to verify that
$$\f d{dt}f(t)=4\times\f{47(1-t^2)^6+616(1-t^2)^3+128}{(8-(1-t^2)^3)^3}.$$
Thus
\begin{align*}&\sum_{k=0}^\infty\f{(63k^2+18k+2)8^k}{(6k+1)\bi{6k}{3k}}
= 2(f(1)-f(0))= 2f(1)=4+\f 32\pi
\end{align*}
as desired. This ends the proof. \qed

\begin{lemma} For any $n\in\N$, we have
\begin{equation}\label{6k+15}
\begin{aligned}&\ \sum_{k=0}^n\f{x^k(9(x-64)k^3+18(x-16)k^2+(11x+208)k+2x+40)}{(6k+1)(6k+5)
\binom{6k}{3k}}
\\=&\ 8+\f{(n+1)(3n+1)(3n+2)x^{n+1}}{(6n+1)(6n+5)\binom{6n}{3n}}.
\end{aligned}
\end{equation}
and
\begin{equation}\label{2k+1}
\begin{aligned}&\ \sum_{k=0}^n\f{x^k(9(x-64)k^3+18(x-48)k^2+(11x-368)k+2x-40)}{(2k+1)(6k+1)(6k+5)
\binom{6k}{3k}}
\\=&\ -8+\f{(n+1)(3n+1)(3n+2)x^{n+1}}{(2n+1)(6n+1)(6n+5)\binom{6n}{3n}}.
\end{aligned}
\end{equation}
\end{lemma}

This can be easily proved by induction on $n$.

\medskip\noindent
{\tt Proof of Theorem \ref{Th-8^k}}.  Putting $x=8$ in \eqref{6k+15} and \eqref{2k+1}, we get the identities
$$\sum_{k=0}^n\f{8^k(63k^3+18k^2-37k-7)}{(6k+1)(6k+5)\binom{6k}{3k}}
=-1-\f{8^n(n+1)(3n+1)(3n+2)}{(6n+1)(6n+5)\binom{6n}{3n}}$$
and
$$\sum_{k=0}^n\f{8^k(63k^3+90k^2+35k+3)}{(2k+1)(6k+1)(6k+5)\binom{6k}{3k}}
=1-\f{8^n(n+1)(3n+1)(3n+2)}{(2n+1)(6n+1)(6n+5)\binom{6n}{3n}}.$$
Letting $n\to+\infty$ we obtain the identities
\begin{equation}\label{-1}\sum_{k=0}^\infty\f{8^k(63k^3+18k^2-37k-7)}{(6k+1)(6k+5)\binom{6k}{3k}}=-1
\end{equation}
and
\begin{equation}\label{1}\sum_{k=0}^\infty\f{8^k(63k^3+90k^2+35k+3)}{(2k+1)(6k+1)(6k+5)\binom{6k}{3k}}=1.
\end{equation}

Since
$$(6k+5)(63k^2+18k+2)-6(63k^3+18k^2-37k-7)=315k^2+324k+52,$$
in view of \eqref{3pi/2} and \eqref{-1} we have
\begin{align*}&\ \sum_{k=0}^\infty\f{8^k(315k^2+324k+52)}{(6k+1)(6k+5)\binom{6k}{3k}}
\\=&\ \sum_{k=0}^\infty\f{8^k(63k^2+18k+2)}{(6k+1)\binom{6k}{3k}}
-6\sum_{k=0}^\infty\f{8^k(63k^3+18k^2-37k-7)}{(6k+1)(6k+5)\binom{6k}{3k}}
\\=&\ \f32\pi+4 -6\times(-1)=\f32\pi+10.
\end{align*}
As
$$(2k+1)(315k^2+324k+52)-10(63k^3+90k^2+35k+3)=63k^2+78k+22,$$
with the aid of \eqref{1} we deduce further that
\begin{align*}&\ \sum_{k=0}^\infty\f{8^k(63k^2+78k+22)}{(2k+1)(6k+1)(6k+5)\binom{6k}{3k}}
\\=&\ \sum_{k=0}^\infty\f{8^k(315k^2+324k+52)}{(6k+1)(6k+5)\binom{6k}{3k}}
-10\sum_{k=0}^\infty\f{8^k(63k^3+90k^2+35k+3)}{(2k+1)(6k+1)(6k+5)\binom{6k}{3k}}
\\=&\ \f32\pi+10-10\times1=\f32\pi.
\end{align*}
This proves the desired \eqref{pi}. \qed

\begin{lemma} For any $n\in\N$, we have
\begin{equation}\label{4096-1}
\sum_{k=0}^n\f{(4536k^3-4500k^2+978k+5)\bi{6k}{3k}}{(2k-1)(6k-1)(6k-5)4096^k}=-\f{\bi{6n}{3n}}{4096^n}.
\end{equation}
and
\begin{equation}\label{4096-2}
\sum_{k=0}^n\f{(4536k^3-4644k^2+1074k+25)\bi{6k}{3k}}{(6k-5)4096^k}=-\f{(2n+1)(6n+5)\bi{6n}{3n}}{4096^n}.
\end{equation}
\end{lemma}
\Proof. It is easy to show \eqref{4096-1} and \eqref{4096-2} by induction on $n$. \qed

\medskip
\noindent{\tt Proof of Theorem \ref{Th1.2}}. Letting $n\to+\infty$ in \eqref{4096-1}, we get
\begin{equation}\label{4500}\sum_{k=0}^\infty\f{(4536k^3-4500k^2+978k+5)\bi{6k}{3k}}{(2k-1)(6k-1)(6k-5)4096^k}=0.
\end{equation}
Since
$$4536k^3-4500k^2+978k+5-18k(252k^2-264k+61)=252k^2-120k+5,$$
combining \eqref{4500} and \eqref{CZ4096} we obtain
$$\sum_{k=0}^\infty\f{(252k^2-120k+5)\bi{6k}{3k}}{(2k-1)(6k-1)(6k-5)}
=0-18\times\f1{12\sqrt3}=-\f{\sqrt3}2.$$
This proves \eqref{252}.

Observe that
\begin{align*}&\ 4536k^3-4500k^2+978k+5
\\=&\ 9(2k-1)(6k-1)(42k+5)-8(252k^2-120k+5).
\end{align*}
Thus, in view of \eqref{4500} and \eqref{252}, we have
\begin{align*}\sum_{k=0}^\infty\f{(42k+5)\bi{6k}{3k}}{(6k-5)4096^k}
=&\ \f 89\sum_{k=0}^\infty\f{(252k^2-120k+5)\bi{6k}{3k}}{(2k-1)(6k-1)(6k-5)4096^k}
\\=&\ \f 89\times\l(-\f{\sqrt3}2\r)=-\f49\sqrt3.
\end{align*}
This proves \eqref{42k+5}.

Letting $n\to+\infty$ in \eqref{4096-2}, we get the identity
\begin{equation}\label{I(6k-5)}\sum_{k=0}^\infty\f{(4536k^3-4644k^2+1074k+25)\bi{6k}{3k}}{(6k-5)4096^k}=0.
\end{equation}
Since
$$4536k^3-4644k^2+1074k+25-3(6k-5)(6k-1)(42k-1)=8(42k+5),$$
combining \eqref{I(6k-5)} with \eqref{42k+5}, we obtain
$$\sum_{k=0}^\infty\f{(6k-1)(42k-1)\bi{6k}{3k}}{4096^k}
=-\f83\sum_{k=0}^\infty\f{(42k+5)\bi{6k}{3k}}{6k-5}=\f{32}{27}\sqrt3.$$
This proves \eqref{42k-1}.

In view of the above, we have completed the proof of Theorem \ref{Th1.2}. \qed

\section{Conjectural series with summands involving $\bi{6k}{3k}$}
\setcounter{equation}{0}
 \setcounter{conjecture}{0}
 \setcounter{theorem}{0}
 \setcounter{proposition}{0}

 \begin{conjecture} We have
 \begin{equation}\sum_{k=0}^\infty\f{8^k((63k^2+78k+22)H_{3k}+12k+8)}{(2k+1)(6k+1)(6k+5)\binom{6k}{3k}}
 =6G-\f32\pi\log2.
 \end{equation}
 \end{conjecture}
 \begin{remark} This is motivated by \eqref{pi}.
 \end{remark}

 \begin{conjecture} We have
 \begin{equation}\label{16^k}\sum_{k=0}^\infty\f{(3k+2)^216^k}{(2k+1)(6k+1)(6k+5)\bi{6k}{3k}}=\f{\pi}{2\sqrt3}.
 \end{equation}
 Also,
 \begin{equation}\sum_{k=0}^\infty\f{(3k+2)^216^kH_k}{(2k+1)(6k+1)(6k+5)\bi{6k}{3k}}=\f 52K-\f{\pi\log432}{6\sqrt3}
 \end{equation}
 and
 \begin{equation}\sum_{k=0}^\infty\f{16^k((3k+2)^2(3H_{3k}-2H_{2k})+\f{2k(3k+2)}{2k+1})}
 {(2k+1)(6k+1)(6k+5)\bi{6k}{3k}}=\f 54K+\f{\pi}{6\sqrt3}\log\f{27}{256}.
 \end{equation}
 \end{conjecture}
 \begin{remark} In the spirit of the proof of Theorem \ref{Th-8^k}, we can show that
 \eqref{16^k} is equivalent to the identity
 $$\sum_{k=0}^\infty\f{16^k(54k^2-36k-12)}{(6k+1)\bi{6k}{3k}}=-4-\f{\pi}{\sqrt3},$$
 which we cannot prove via integration.
 \end{remark}

 \begin{conjecture} We have
 \begin{equation}\sum_{k=0}^\infty\f{(21k^2+27k+8)(H_{2k}+2H_k)+\f{(3k+1)(5k+4)}{6k+3}}{(2k+1)(6k+1)(6k+5)\bi{6k}{3k}}
 =14K-\f{16\log3}{3\sqrt3}\pi
 \end{equation}
 and
\begin{equation}\sum_{k=0}^\infty\f{(21k^2+27k+8)H_{3k}+2k+13/9}{(2k+1)(6k+1)(6k+5)\bi{6k}{3k}}
=\f{8}3K-\f {8\log3}{9\sqrt3}\pi.
\end{equation}
 \end{conjecture}
 \begin{remark} This is motivated by the identity \eqref{CZ115}.
 \end{remark}

  \begin{conjecture} We have
 \begin{equation}\sum_{k=0}^\infty\f{27^k((111k^2+133k+36)(H_{2k}-3H_k)-\f{(3k+1)(3k+4)}{2k+1})}{(2k+1)(6k+1)(6k+5)\bi{6k}{3k}}
 =\f{16\log3}{\sqrt3}\pi-48K
 \end{equation}
 and
\begin{equation}\sum_{k=0}^\infty\f{27^k((111k^2+133k+36)H_{3k}+30k+19}{(2k+1)(6k+1)(6k+5)\bi{6k}{3k}}
=12K.
\end{equation}
 \end{conjecture}
 \begin{remark} This is motivated by the identity
 $$\sum_{k=0}^\infty\f{27^k(111k^2+133k+36)}{(2k+1)(6k+1)(6k+5)\bi{6k}{3k}}=\f{16\pi}{3\sqrt3}$$
 given by \cite[Example 58]{CZ}.
 \end{remark}

 \begin{conjecture} We have
 \begin{equation}\begin{aligned}&\ \sum_{k=0}^\infty\f{((6k-1)(42k-1)(3H_{3k}-10H_{2k}+4H_k)-1836k-50)\bi{6k}{3k}}{4096^k}
\\&\qquad= -\f{32\sqrt3}{27}\l(32+14\log2-9\log3\r)
 \end{aligned}
 \end{equation}
 and
 \begin{equation}\begin{aligned}&\ \sum_{k=0}^\infty\f{((6k-1)(42k-1)(2H_{6k}-H_{3k})+24k+8/3)\bi{6k}{3k}}{4096^k}
\\&\qquad= \f{32\sqrt3}{81}\l(8+\log\f{64}{27}\r).
 \end{aligned}
 \end{equation}
 \end{conjecture}
 \begin{remark} This is motivated by the identity \eqref{42k-1}.
 \end{remark}

 \begin{conjecture} We have
 \begin{equation}\label{216}\sum_{k=0}^\infty\f{(57k^2+2k+1)\bi{6k}{3k}}{216^k}=9\sqrt3
 \end{equation}
 and
 \begin{equation}\sum_{k=0}^\infty\f{((57k^2+2k+1)(2H_{6k}-H_{3k})+14k+8/3)\bi{6k}{3k}}{216^k}=
 3\sqrt3(8+3\log3).
 \end{equation}
 \end{conjecture}
 \begin{remark} In the spirit of the proof of Theorem \ref{Th1.2}, \eqref{216} has the following equivalent form:
 \begin{equation}\sum_{k=0}^\infty\f{(57k^2-40k+5)\bi{6k}{3k}}{(2k-1)(6k-1)(6k-5)216^k}=-\f1{\sqrt3}.
 \end{equation}
 \end{remark}

 \begin{conjecture} We have
 \begin{equation}\label{256}\sum_{k=0}^\infty\f{(6k-1)\bi{6k}{3k}}{256^k}=-\f{2}{3\sqrt3}
 \end{equation}
 and
 \begin{equation}\sum_{k=0}^\infty\f{((6k-1)(3H_{3k}-H_k)-4)\bi{6k}{3k}}{256^k}
 =\f{4\sqrt3}9(2\log2-3).
 \end{equation}
 Also,
 \begin{equation}\sum_{k=0}^\infty\f{k\bi{6k}{3k}(2H_{2k}-H_k-\f{2}{6k-3})}{(6k-1)(6k-5)256^k}
 =\f{\sqrt3}{324}\l(6+\log\f4{27}\r)
 \end{equation}
 and
 \begin{equation}\sum_{k=0}^\infty\f{k\bi{6k}{3k}(3H_{3k}-H_k-4/(6k-3))}{(6k-1)(6k-5)256^k}
 =\f{\sqrt3}{54}(3-2\log2).
 \end{equation}
 \end{conjecture}
 \begin{remark} In the spirit of the proof of Theorem \ref{Th1.2}, \eqref{256} has the following equivalent form:
 \begin{equation}\sum_{k=0}^\infty\f{k\bi{6k}{3k}}{(6k-1)(6k-5)256^k}=\f{\sqrt3}{108}.
 \end{equation}
 \end{remark}

 \begin{conjecture} We have
 \begin{equation}\label{512}\sum_{k=0}^\infty\f{(84k^2-4k+1)\bi{6k}{3k}}{512^k}=4\sqrt2
 \end{equation}
 and
 \begin{equation}\sum_{k=0}^\infty\f{((84k^2-4k+1)(2H_{6k}-H_{3k})+16k+8/3)\bi{6k}{3k}}{512^k}
 =\f{4\sqrt2}3(8+3\log2).
 \end{equation}
 Also,
 \begin{equation}\sum_{k=0}^\infty\f{\bi{6k}{3k}((42k-5)(H_{6k}-\f{H_{3k}}2)-\f{1224k^3-1308k^2+366k-25}{(2k-1)(6k-1)(6k-5)})}
 {(6k-1)(6k-5)512^k}=-\f{2+\log2}{2\sqrt2}.
 \end{equation}
 \end{conjecture}
 \begin{remark} In the spirit of the proof of Theorem \ref{Th1.2}, \eqref{512} has the following equivalent form:
  \begin{equation}\sum_{k=0}^\infty\f{(42k-5)\bi{6k}{3k}}{(6k-1)(6k-5)512^k}=-\f{\sqrt2}2.
  \end{equation}
 \end{remark}

 \begin{conjecture} We have
 \begin{equation}\label{-512}\sum_{k=0}^\infty\f{(6k-1)^2\bi{6k}{3k}}{(-512)^k}=\f4{27}\sqrt6,
 \end{equation}
 \begin{equation}\sum_{k=0}^\infty\f{(6k-1)^2(3H_{3k}-4H_{2k}+H_k)\bi{6k}{3k}}{(-512)^k}
 =\f{4\sqrt6}9\log\f34,
 \end{equation}
 and
 \begin{equation}\sum_{k=0}^\infty\f{((6k-1)^2(2H_{6k}-H_{3k})+8/3)\bi{6k}{3k}}{(-512)^k}
 =\f{4\sqrt6}{81}\l(8+3\log\f23\r).
 \end{equation}
 Also,
 \begin{equation}\begin{aligned}&\ \sum_{k=0}^\infty\f{\bi{6k}{3k}((36k^2-12k+5)H(k)+\f{16P(k)}{(6k-1)(6k-3)(6k-5)})}
 {(2k-1)(6k-1)(6k-5)(-512)^k}
\\&\qquad= \sqrt6\l(2\log2-\f{\log3}2-8\r),
 \end{aligned}\end{equation}
 where $H(k)=4H_{6k}+H_{3k}-4H_{2k}+H_k$ and $P(k)=936k^3-1476k^2+622k-75$.
 \end{conjecture}
\begin{remark}  In the spirit of the proof of Theorem \ref{Th1.2}, \eqref{-512} has the following equivalent form:
 \begin{equation}
 \sum_{k=0}^\infty\f{(36k^2-12k+5)\bi{6k}{3k}}{(2k-1)(6k-1)(6k-5)(-512)^k}=-\f{\sqrt6}2.
 \end{equation}
 \end{remark}

\section{Duality principle for irrational series}
 \setcounter{equation}{0}
 \setcounter{conjecture}{0}
 \setcounter{theorem}{0}
 \setcounter{proposition}{0}

 Ramanujan's list (cf. \cite{R}) of series for $1/\pi$ contains only one irrational series:
 $$\sum_{k=0}^\infty\l(k+\f{31}{270+48\sqrt5}\r)\f{\bi{2k}k^3}{(2^4(\sqrt5+1)^8)^k}=\f{16}{3(5+7\sqrt5)\pi},$$
 i.e.,
 \begin{equation}\label{R1}\sum_{k=0}^\infty\l(6k(7\sqrt5+5)+5\sqrt5-1\r)\f{\bi{2k}k^3}{(12+4\sqrt5)^{4k}}
 =\f{32}{\pi}.\end{equation}
 In 1988 J. M. Borwein and P. B. Borwein \cite{BB} gave 7 pairs of new irrational series of Ramanujan's type; for example,
 \begin{equation}\label{19}\sum_{k=0}^\infty\f{\bi{2k}k^2\bi{4k}{2k}(40k(5457\sqrt{2}+693121)-869892\sqrt{2}+1877581)}{19^{2k}(12(481+340\sqrt{2}))^{4k}}
 =466578\f{\sqrt{19}}{\pi}
 \end{equation}
 and
 \begin{equation*}\sum_{k=0}^\infty\f{\bi{2k}k^2\bi{4k}{2k}(40k(693121-5457\sqrt{2})+869892\sqrt{2}+1877581)}
  {19^{2k}(12(481-340\sqrt{2}))^{4k}}
 =2332890\f{\sqrt{19}}{\pi}.
 \end{equation*}
 Note that the series in \eqref{19} converges quickly with the geometric rate about $3.996\times10^{-17}$.

 In 2012 A. M. Aldawoud \cite{A} did a systematic search via Maple, and guessed 76 new irrational series of the Ramanujan type. However, some of her tables (such as the fourth row of Table 3.4)
 contain false data which might be caused by typos; for example, the tenth row of \cite[Table 3.1]{A}
 indicates that for
 $$x=\f{8145698488}{79267303539275}-\f{570145303\sqrt7}{14679130284125}
 \ \t{and} \ \lambda=\f{90352}{1274007}-\f{59000\sqrt7}{3822021}$$
 we should have
 $$\pi\sqrt{7(1-1728x)}\sum_{k=0}^\infty\bi{2k}k\bi{3k}k\bi{6k}{3k}x^k(k+\lambda)=\f14,$$
 but our calculation shows that
 $$\pi\sqrt{7(1-1728x)}\sum_{k=0}^\infty\bi{2k}k\bi{3k}k\bi{6k}{3k}x^k(k+\lambda)
 \approx  0.24999999999474.$$

 In 2023 J. M. Campbell \cite{Cam} proved some irrational Ramanujan-type series of level $3$ (with summands involving $\bi{2k}k^2\bi{3k}k$) conjectured by Aldawoud \cite{A}.

 Let $d$ be any positive integer which is not a square.
 By Galois theory,  the Galois group of the field extension $\Q(\sqrt d)/\Q$ has exactly two elements, and the nontrivial Galois automorphism
 $\sigma_d$ in $\mathrm{Gal}(\Q(\sqrt d)/\Q)$ sends $a+b\sqrt d$ (with $a,b\in\Q$) to its conjugate $a-b\sqrt d$. Inspired by this, for an irrational series
 $\sum_{k=0}^\infty a_k$ with $a_k\in\Q(\sqrt d)$, we define its {\it dual} by
$$\lim_{n\to+\infty}\sigma_d\(\sum_{k=0}^n a_k\)=\sum_{k=0}^\infty\sigma_d(a_k)$$
if the series $\sum_{k=0}^\infty\sigma_d(a_k)$ converges.

The dual series of the Ramanujan series in \eqref{R1}
 is
 $$\sum_{k=0}^\infty\l(6k(-7\sqrt5+5)-5\sqrt5-1\r)\f{\bi{2k}k^3}{(12-4\sqrt5)^{4k}},$$
and the identity \eqref{R1} has the following dual:
 \begin{equation}\label{R2}\sum_{k=0}^\infty\l(6k(7\sqrt5-5)+5\sqrt5+1\r)\f{\bi{2k}k^3}{(12-4\sqrt5)^{4k}}
 =\f{96}{\pi},\end{equation}
 which was noted in \cite{G-AMM}.

 In the spirit of the author's recent papers \cite{harmonic,SSM,S24,SII}, we pose the following conjecture motivated by \eqref{R1} and \eqref{R2}.

 \begin{conjecture} \label{Conj1.1} We have
 \begin{equation}\sum_{k=0}^\infty\l(6k(5+7\sqrt5)+5\sqrt5-1)H(k)
 +\f{60(7-3\sqrt5)}{2k+1}\r)\f{\bi{2k}k^3}{(12+4\sqrt5)^{4k}}=\f{16}3\pi
 \end{equation}
 and
 \begin{equation}\sum_{k=0}^\infty\l(6k(5-7\sqrt5)-5\sqrt5-1) H(k)
 +\f{60(7+3\sqrt5)}{2k+1}\r)\f{\bi{2k}k^3}{(12-4\sqrt5)^{4k}}=1136\pi.
 \end{equation}
 where
 $$ H(k):=35H_k^{(2)}-136H_{2k}^{(2)}.$$
 \end{conjecture}

 For any integer $d\eq0,1\pmod4$, we adopt the notation
$$L_{d}(2):=L\l(2,\l(\f{d}{\cdot}\r)\r)=\sum_{k=1}^\infty\f{(\f dk)}{k^2}$$
with $(\f dk)$ the Kronecker symbol. Thus $G=L_{-4}(2)$ and $K=L_{-3}(2)$.

For convenience, let us call the series of the form
$$\sum_{k=1}^\infty\f{(ak-b)m^k}{k^3c_k}=r_1L_{n'}(2)+r_2\sqrt{d}L_{(dn)'}(2)$$
with $am$ nonzero and $c_k$ among \eqref{three} a {\it Zeilberger-type series}, where
$d\in\Z^+=\{1,2,3,\ldots\}$, $n\in\Z$, $a,b,m\in\Q(\sqrt d)$, $r_1,r_2\in\Q$, and $c'$ with $c\in\Z$ is defined by
$$c'=\begin{cases}c&\t{if}\ c\eq0,1\pmod 4,\\4c&\t{otherwise}.\end{cases}$$
When $d=1$, we call such series {\it rational series of the Zeilberger type}.
It seems that \eqref{Z}-\eqref{S2} are the only rational series of Zeilberger's type. Table 5 of \cite{GR} provides ten irrational series of Zeilberger's type one of which is the identity
\begin{equation}\label{GR5}
\sum_{k=1}^\infty\f{3(16\sqrt5+35)k-4(5\sqrt5+11)}{k^3\bi{2k}k^3}\l(\f{1-\sqrt5}2\r)^{8k}=\f{\pi^2}{30}.
\end{equation}
We find the dual of this identity:
\begin{equation}\label{GR-5}
\sum_{k=1}^\infty\f{3(16\sqrt5-35)k-4(5\sqrt5-11)}{k^3\bi{2k}k^3}\l(\f{1+\sqrt5}2\r)^{8k}
=\f{71}{30}\pi^2,
\end{equation}
which will be proved in a forthcoming paper joint with Yajun Zhou.

We make the following conjecture similar to Conjecture \ref{Conj1.1}.

\begin{conjecture} We have
\begin{equation}\sum_{k=1}^\infty\f{((1-\sqrt5)/2)^{8k}}{k^3\bi{2k}k^3}
\l(\l(3(16\sqrt5+35)k-4(5\sqrt5+11)\r)\mathcal H(k)+\f{15}k\r)=\f{\pi^4}{72}
\end{equation}
and
\begin{equation}\sum_{k=1}^\infty\f{((1+\sqrt5)/2)^{8k}}{k^3\bi{2k}k^3}
\l(\l(3(35-16\sqrt5)k+4(5\sqrt5-11)\r)\mathcal H(k)+\f{15}k\r)=\f{433}{72}\pi^4.
\end{equation}
where
$$\mathcal H(k):=4H_{2k-1}^{(2)}-35H_{k-1}^{(2)}.$$
\end{conjecture}

Let $d$ be a positive integer which is not a square. Let $c_k\ (k\in\N)$
be one of the four products in \eqref{three}, and set $c=\lim_{k\to\infty}\root k\of{c_k}$.
For a Ramanujan-type series
$$\sum_{k=0}^\infty\f{(ak+b)c_k}{m^k}=\f{r\sqrt n}{\pi}$$
with $n$  a positive algebraic integer in $\Q(\sqrt d)$, $a,b,m\in\Q(\sqrt d)$, $r\in\Q$ and $amr\not=0$,
if $|\sigma_d(m)|<c$ then
we define the  {\it dual} of the Ramanujan series $\sum_{k=0}^\infty\f{(ak+b)c_k}{m^k}$
as the Zeilberger-type series
$$\sum_{k=1}^\infty\sigma_d\l(\f{(ak-b)m^k}{k^3c_k}\r).$$

\begin{conjecture}[Duality Principle] Let $d$ be a positive integer which is not a square. And let $c_k\ (k\in\N)$
be one of the four products in \eqref{three}. Let $a,b,m\in\Q(\sqrt d)$ with $a\not=0$ and $m\not\in\Q$. Let $c$ denote the positive integer $\lim_{k\to+\infty}\root k\of{c_k}$.

{\rm (i)} Suppose that we have a Ramanujan-type series
$$\sum_{k=0}^\infty\f{(ak+b)c_k}{m^k}=\f{r\sqrt n}{\pi}$$
where $r\in\Q\sm\{0\}$, and $n$ is a positive algebraic integer in $\Q(\sqrt d)$.
If $\sum_{k=0}^\infty\sigma_d\l(\f{(ak+b)c_k}{m^k}\r)$ converges (i.e., $|\sigma_d(m)|>c$)
and $\sigma_d(n)>0$, then
$$\sum_{k=0}^\infty\sigma_d\l(\f{(ak+b)c_k}{m^k}\r)=\f{r^*\sqrt {\sigma_d(n)}}{\pi}$$
for some rational number $r^*$.
If $\sum_{k=0}^\infty\sigma_d(\f{(ak+b)c_k}{m^k})$ diverges (i.e. $|\sigma_d(m)|<c$),
and $m<0$ and $n\in\Z^+$, then
we have a Zeilberger-type series
$$\sum_{k=1}^\infty\sigma_d\l(\f{(ak-b)m^k}{k^3c_k}\r)=r_1L_{(-n)'}(2)+r_2\sqrt{d}L_{(-dn)'}(2)$$
for some $r_1,r_2\in\Q$.

{\rm (ii)} Suppose that we have a Zeilberger-type series
$$\sum_{k=1}^\infty\f{(ak-b)m^k}{k^3c_k}=r_1L_{n'}(2)+r_2\sqrt{d}L_{(dn)'}(2),$$
where $n\in\Z$ and $r_1,r_2\in\Q$. If the series $\sum_{k=1}^\infty\sigma_d(\f{(ak-b)m^k}{k^3c_k})$ converges (i.e., $|\sigma_d(m)|<c$), then
$$\sum_{k=1}^\infty\sigma_d\l(\f{(ak-b)m^k}{k^3c_k}\r)=r_1^*L_{n'}(2)+r_2^*\sqrt{d}L_{(dn)'}(2)$$
for some $r_1^*,r_2^*\in\Q$.
\end{conjecture}

\section{New irrational series for $1/\pi$}
 \setcounter{equation}{0}
 \setcounter{conjecture}{0}
 \setcounter{theorem}{0}
 \setcounter{proposition}{0}

 In this section, we apply the Duality Principle to find 16 new irrational series for $1/\pi$.

 \begin{conjecture}\label{Conj2.1} We have
 \begin{equation}\label{375}\sum_{k=0}^\infty\f{\bi{2k}k^2\bi{3k}k(9k(51-11\sqrt6)+2(54-19\sqrt6))}
 {(27(37102+15147\sqrt6))^k}
 =\f{375}{8\pi},
 \end{equation}
 \begin{equation}\label{250}\sum_{k=0}^\infty\f{\bi{2k}k^2\bi{3k}k(2k(182-27\sqrt6)+3(16-\sqrt6))}
 {8^k(13-7\sqrt6)^{2k}}
 =\f{250\sqrt2}{\pi},
 \end{equation}
 \begin{equation}\label{162}\sum_{k=0}^\infty\f{\bi{2k}k^2\bi{3k}k(30k(7\sqrt{10}-2)+25\sqrt{10}+16)}
 {108^k(14-5\sqrt{10})^{2k}}
 =\f{162\sqrt6}{\pi},
 \end{equation}
 \begin{equation}\label{slow}\sum_{k=0}^\infty\f{\bi{2k}k^2\bi{3k}k(15k(377\sqrt{3}-378)
 +4(122\sqrt{3}-27))}{(114-63\sqrt{3})^{3k}}
 =\f{11979}{\pi},
 \end{equation}
 \begin{equation}\label{34}\sum_{k=0}^\infty\f{\bi{2k}k^2\bi{3k}k(30k(59228-7\sqrt{34})+125792+7157\sqrt{34})}{108^k(1876-325\sqrt{34})^{2k}}
 =215622\f{\sqrt6}{\pi},
 \end{equation}
 and
 \begin{equation}\label{265}\sum_{k=0}^\infty\f{\bi{2k}k^2\bi{3k}k
 (3k(20801\sqrt{265}-24115)+7049\sqrt{265}-82795)}{(-15)^k(60(341\sqrt{265}+5551))^{2k}}
 =25920\f{\sqrt{15}}{\pi}.
 \end{equation}
 \end{conjecture}
 \begin{remark} \eqref{250}--\eqref{34} are duals of the conjectural identities
 $$\sum_{k=0}^\infty\f{\bi{2k}k^2\bi{3k}k(2k(182+27\sqrt6)+3(16+\sqrt6))}{8^k(13+7\sqrt6)^{2k}}
 =\f{125\sqrt2}{\pi},$$
 $$ \sum_{k=0}^\infty\f{\bi{2k}k^2\bi{3k}k(30k(7\sqrt{10}+2)+25\sqrt{10}-16)}{108^k(14+5\sqrt{10})^{2k}}
 =\f{81\sqrt6}{\pi},$$
 $$ \sum_{k=0}^\infty\f{\bi{2k}k^2\bi{3k}k(15k(377\sqrt{3}+378)+4(122\sqrt{3}+27))}
 {(114+63\sqrt{3})^{3k}}
 =\f{11979}{4\pi},$$
 and
 $$\sum_{k=0}^\infty\f{\bi{2k}k^2\bi{3k}k(30k(59228+7\sqrt{34})+125792-7157\sqrt{34})}{108^k(1876+325\sqrt{34})^{2k}}
 =107811\f{\sqrt6}{\pi}$$
 given by rows two, five, nine and eleven of \cite[Table 3.6]{A}, respectively.
 Also, \eqref{265} is the dual of the conjectural identity
 \begin{equation*}\sum_{k=0}^\infty\f{\bi{2k}k^2\bi{3k}k
 (3k(20801\sqrt{265}+24115)+7049\sqrt{265}+82795)}{(-15)^k(60(341\sqrt{265}-5551))^{2k}}
 =360^2\f{\sqrt{15}}{\pi}
 \end{equation*}
 given by the eighth row of \cite[Table 3.8]{A}.
 \end{remark}

 \begin{conjecture} \label{Conj2.2} We have
 \begin{equation}\label{121}
 \sum_{k=0}^\infty\f{\bi{2k}k^2\bi{4k}{2k}(5k(32\sqrt3-17)+18\sqrt3-2)}{(36(746-425\sqrt3))^k}
 =\f{121\sqrt3}{\pi},
 \end{equation}
 \begin{equation}\label{441}
 \sum_{k=0}^\infty\f{\bi{2k}k^2\bi{4k}{2k}((352\sqrt2+57)k+2(23\sqrt2+15))}{(162(884-627\sqrt2))^k}
 =\f{441}{2\pi},
 \end{equation}
 \begin{equation}\label{7Pi}\sum_{k=0}^\infty\f{\bi{2k}k^2\bi{4k}{2k}2^k(35k(2176-13\sqrt7)+8(98\sqrt7+841))}{(9(325-119\sqrt7))^{2k}}
 =\f{171^2}{\pi},
 \end{equation}
  \begin{equation}\label{98}
 \sum_{k=0}^\infty\f{\bi{2k}k^2\bi{4k}{2k}(40k(3\sqrt2+19)-12\sqrt2+71)}{(48(57+40\sqrt2))^{2k}}
 =\f{98\sqrt3}{\pi},
 \end{equation}
\begin{equation}\label{340}\sum_{k=0}^\infty\f{(-1)^k\bi{2k}k^2\bi{4k}{2k}(20k(323-124\sqrt3)+1329-708\sqrt3)}
{(32(340\sqrt3+589))^{2k}}
=\f{968}{3\pi},
\end{equation}
 \begin{equation}\label{384}\sum_{k=0}^\infty\f{\bi{2k}k^2\bi{4k}{2k}(20k(13\sqrt{33}+11)+33\sqrt{33}+119)}
 {(-3072)^k(75-13\sqrt{33})^{2k}}
 =384\f{\sqrt3}{\pi},
 \end{equation}
and
 \begin{equation}\label{57Pi}\sum_{k=0}^\infty\f{\bi{2k}k^2\bi{4k}{2k}(260k(\sqrt{57}+323)+513\sqrt{57}+7331)}
 {(18432(71825\sqrt{57}-542267))^k}=18816\f{\sqrt3}{\pi}.
 \end{equation}
 \end{conjecture}
 \begin{remark} \eqref{121}, \eqref{441}, \eqref{7Pi} and \eqref{98}
 are the duals of the conjectural identities
 $$\sum_{k=0}^\infty\l(5k(17+32\sqrt 3)+2(9\sqrt3+1)\r)\bi{2k}k^2\bi{4k}{2k}\l(\f{746-425\sqrt3}{527076}\r)^k
 =\f1{2\sqrt3\,\pi},$$
 $$\sum_{k=0}^\infty\l((352\sqrt 2-57)k+2(23\sqrt2-15)\r)\bi{2k}k^2\bi{4k}{2k}\l(\f{627\sqrt2-884}{777924}\r)^k
 =\f{441}{4\pi},$$
 $$\sum_{k=0}^\infty\f{\bi{2k}k^2\bi{4k}{2k}2^k(35k(2176-13\sqrt7)+8(98\sqrt7+841)}
 {(9(325+119\sqrt7))^{2k}}
 =\f{171^2}{2\pi},$$
 and
 \begin{equation*}\sum_{k=0}^\infty\f{\bi{2k}k^2\bi{4k}{2k}(40k(19-3\sqrt2)+12\sqrt2+71)}
 {(48(57-40\sqrt2))^{2k}}=\f{294\sqrt3}{\pi},\end{equation*}
 given by rows two, four, five and six of \cite[Table 3.3]{A}, respectively.
 As the six row of \cite[Table 3.3]{A} actually gave a wrong value of $\lambda$, we find the last identity by using the PSLQ algorithm (cf. \cite{PSLQ}).
 Moreover,
  \eqref{384} is the dual of the conjectural identities
  \begin{equation*}\sum_{k=0}^\infty\f{\bi{2k}k^2\bi{4k}{2k}(20k(13\sqrt{33}-11)+33\sqrt{33}-119)}{(-3072)^k(75+13\sqrt{33})^{2k}}
 =128\f{\sqrt3}{\pi},
 \end{equation*}
 given by the third row of \cite[Table 3.5]{A}, and \eqref{57Pi} is the dual of the conjectural identities
 \begin{equation*}\sum_{k=0}^\infty\f{\bi{2k}k^2\bi{4k}{2k}(260k(323-\sqrt{57})-513\sqrt{57}+7331)}
 {(-18432(71825\sqrt{57}+542267))^k}=6272\f{\sqrt3}{\pi}
 \end{equation*}
 given by fifth row of \cite[Table 3.5]{A}.
 Note also that the dual of \eqref{340} is the Zeilberger-type identity
 given by the last row of \cite[Table 5]{GR}.
 \end{remark}

\begin{conjecture}\label{Conj2.3} We have
\begin{equation}\label{32}\begin{aligned}&\sum_{k=0}^\infty\f{\bi{2k}k\bi{3k}k\bi{6k}{3k}(92158\sqrt5 k+5(1393\sqrt5+1600))}
{96^{3k}(165393\sqrt5-369830)^k}\\&\quad =\f{32}{5\pi}\sqrt{30(3777190-165393\sqrt5)},
\end{aligned}
\end{equation}
\begin{equation}\label{195}\begin{aligned}&\sum_{k=0}^\infty\f{\bi{2k}k\bi{3k}k\bi{6k}{3k}(20124k+125\sqrt{13}
+1339)}
{60^{3k}(15965-4428\sqrt{13})^k}
\\&\qquad=\f{10}{3\pi}\sqrt{195(14230-41\sqrt{13})},
\end{aligned}
\end{equation}
and
\begin{equation}\label{135}\begin{aligned}
&\sum_{k=0}^\infty\f{\bi{2k}k\bi{3k}k\bi{6k}{3k}2^k(13629k+220\sqrt5+1052)}
{(135(274207975-122629507\sqrt5))^k}
\\&\qquad=\f3{2\pi}\sqrt{3(5607725+661856\sqrt5)}.
\end{aligned}
\end{equation}
\end{conjecture}
\begin{remark} \eqref{32}, \eqref{195} and \eqref{135} are the duals of the conjectural identities
\begin{align*}&\sum_{k=0}^\infty\f{\bi{2k}k\bi{3k}k\bi{6k}{3k}(92158\sqrt5 k+5(1393\sqrt5-1600))}
{(-96)^{3k}(165393\sqrt5+369830)^k}\\&\quad =\f{32}{15\pi}\sqrt{30(3777190+165393\sqrt5)},
\end{align*}
\begin{align*}&\sum_{k=0}^\infty\f{\bi{2k}k\bi{3k}k\bi{6k}{3k}(20124k-125\sqrt{13}
+1339)}
{60^{3k}(15965+4428\sqrt{13})^k}=\f{5}{3\pi}\sqrt{195(14230+41\sqrt{13})}
\end{align*}
and
\begin{align*}
&\sum_{k=0}^\infty\f{\bi{2k}k\bi{3k}k\bi{6k}{3k}2^k(13629k-220\sqrt5+1052)}
{(135(274207975+122629507\sqrt5))^k}
\\&\qquad=\f1{2\pi}\sqrt{3(5607725-661856\sqrt5)}
\end{align*}
given by the second row of \cite[Tsable 3.2]{A}, rows five and six of \cite[Table 3.1]{A}, respectively.
\end{remark}

\section{New irrational series of Zeilberger's type}
 \setcounter{equation}{0}
 \setcounter{conjecture}{0}
 \setcounter{theorem}{0}
 \setcounter{proposition}{0}

In this section, we apply the Duality Principle to find 10 new irrational series for $1/\pi$.
They provide fast converging series for computing the $L$-values
$$L_{-8}(2),\ L_{-11}(2),\ L_{-24}(2),\ L_{-39}(2),\ L_{-68}(2)\ \t{and}\ L_{-111}(2).$$

\begin{conjecture}\label{Conj3.1} We have
\begin{equation}\sum_{k=1}^\infty\f{(24(12-7\sqrt3))^k}{k^3\bi{2k}k^2\bi{3k}k}
\l((3\sqrt3+7)k-\sqrt3-2\r)=\f56\l(9\sqrt3 K-16G\r),
\end{equation}
\begin{equation}\begin{aligned}&\sum_{k=1}^\infty\f{(54(265-153\sqrt3))^k}{k^3\bi{2k}k^2\bi{3k}k}
\l(3k(17\sqrt3+27)-16\sqrt3-27\r)
\\&\qquad=135\l(G-\f{11}{16}\sqrt3 K\r),
\end{aligned}\end{equation}
\begin{equation}\label{-24}\begin{aligned}&\sum_{k=1}^\infty\f{(27(37102-15147\sqrt6))^k}{k^3\bi{2k}k^2\bi{3k}k}
\l(9k(11\sqrt6+51)-2(19\sqrt6+54)\r)
\\&\qquad=\f{5625}2\l(\sqrt6 L_{-24}(2)-3G\r),
\end{aligned}
\end{equation}
\begin{equation}\begin{aligned}&\sum_{k=1}^\infty\f{(-1728)^k(18-5\sqrt{13})^{2k}}{k^3\bi{2k}k^2\bi{3k}k}
\l(15k(9\sqrt{13}+26)-39\sqrt{13}-134)\r)
\\&\qquad=54\l(80K-13\sqrt{13}L_{-39}(2)\r),
\end{aligned}
\end{equation}
\begin{equation}\label{320}\begin{aligned}&\sum_{k=1}^\infty\f{(32(91\sqrt{33}-523))^{k}}{k^3\bi{2k}k^2\bi{3k}k}
\l((91\sqrt{33}+891)k-33\sqrt{33}-225\r)
\\&\qquad=320\l(\f{11}3\sqrt{33}L_{-11}(2)-27K\r),
\end{aligned}
\end{equation}
\begin{equation}\label{111}\begin{aligned}&\sum_{k=1}^\infty\f{(-3)^k(24(2737-450\sqrt{37}))^{2k}(15k(2898\sqrt{37}+3145)-6438\sqrt{37}-30355)}
{k^3\bi{2k}k^2\bi{3k}k}
\\&\qquad=71874\l(380K-37\sqrt{37}L_{-111}(2)\r),
\end{aligned}
\end{equation}
and
\begin{equation}\label{145}\begin{aligned}&\sum_{k=1}^\infty\f{(-3/5)^k(12m)^{2k}}{k^3\bi{2k}k^2\bi{3k}k}
\l(15k(323\sqrt{145}+13195)-2407\sqrt{145}-35375\r)
\\&\qquad=23328\l(29\sqrt{145}L_{-87}(2)-375L_{-15}(2)\r)
\end{aligned}
\end{equation}
where $m=2975-247\sqrt{145}$.
\end{conjecture}
\begin{remark} \eqref{320} provides a fast way to compute $L_{-11}(2)$.
\eqref{-24} is the dual of \eqref{375},  and
\eqref{145} is the dual of the conjectural identity
$$\sum_{k=0}^\infty\f{\bi{2k}k^2\bi{3k}k(15k(13195-323\sqrt{145})+35375-2407\sqrt{145})}{(-3/5)^k(12(2975+247\sqrt{145}))^{2k}}
=\f{5184\sqrt{15}}{\pi}$$
which corresponds to the corrected form of the fourth row of \cite[Table 3.8]{A} containing typos.
\end{remark}

\begin{conjecture}\label{Conj3.3} We have
\begin{equation}\label{-8}\begin{aligned}&\sum_{k=1}^\infty\f{(8(457-325\sqrt2))^k}{k^3\bi{2k}k^2\bi{4k}{2k}}
\l(5k(13\sqrt2+32)-3(6\sqrt2+11)\r)
\\&\qquad\ \ =490\sqrt2L_{-8}(2)-980G
\end{aligned}
\end{equation}
and
\begin{equation}\label{11}\begin{aligned}&\sum_{k=1}^\infty\f{(-4)^k(84(1121-338\sqrt{11}))^{4k}
(5980ak-b)}{k^3\bi{2k}k^2\bi{4k}{2k}}
\\&\qquad\  =36120^2\l(36G-11\sqrt{11}L_{-11}(2)\r),
\end{aligned}
\end{equation}
where
$$a=11092\sqrt{11}+19437\ \ \t{and}\ \ b=11937508\sqrt{11}+37515813.$$
Also,
\begin{equation}\label{17Z}\begin{aligned}&\sum_{k=1}^\infty\f{(-4)^k(12(41-10\sqrt{17}))^{4k}}{k^3\bi{2k}k^2\bi{4k}{2k}}
\l(140k(92\sqrt{17}+2091)-6868\sqrt{17}-36591\r)
\\&\qquad\ \ =415812(17\sqrt{17}L_{-68}(2)-90G).
\end{aligned}
\end{equation}
\end{conjecture}
\begin{remark} \eqref{-8} is the dual of the conjectural identity
$$\sum_{k=0}^\infty\f{\bi{2k}k^2\bi{4k}{2k}}{(8(457+325\sqrt2))^k}
\l(5k(32-13\sqrt2)+3(11-6\sqrt2)\r)=\f{49}{2\pi}$$
given by the first row of \cite[Table 3.3]{A}, and
\eqref{11} is the dual of the Ramanujan-type identity
$$\sum_{k=0}^\infty\f{\bi{2k}k^2\bi{4k}{2k}(5980a'k+b')}{(-4)^k((84(1121+338\sqrt{11}))^{4k}}
=1806^2\f2{\pi}$$
with $a'=11092\sqrt{11}-19437$ and $b'=11937508\sqrt{11}-37515813$ listed in the last row of \cite[Table 3.5]{A} and first proved in \cite{BB}. Also, \eqref{17Z} is the dual of the conjectural identity
$$\sum_{k=0}^\infty\f{\bi{2k}k^2\bi{4k}{2k}(140k(2091-92\sqrt{17})+36591-6868\sqrt{17})}
{(-4)^k(12(41+10\sqrt{17}))^{4k}}=\f{25992}{\pi}$$
given by the 7th row of \cite[Table 3.5]{A}.
\end{remark}

\end{document}